\input epsf
\documentstyle{amsppt}
\pagewidth{5.4truein}\hcorrection{0.55in}
\pageheight{7.5truein}\vcorrection{0.75in}
\TagsOnRight
\NoRunningHeads
\catcode`\@=11
\def\logo@{}
\footline={\ifnum\pageno>1 \hfil\folio\hfil\else\hfil\fi}
\topmatter
\title The equivalence between enumerating cyclically symmetric, self-complementary\\
and totally symmetric, self-complementary plane partitions
\endtitle
\author Mihai Ciucu\endauthor
\affil
  Institute for Advanced Study\\
  School of Mathematics\\
  Princeton, NJ 08540
\endaffil
\date February, 1998\enddate
\abstract
We prove that the number of cyclically symmetric,
self-complementary plane partitions contained in a cube of side $2n$ equals the square
of the number of totally symmetric, self-complementary plane partitions contained in
the same cube, without explicitly evaluating either of these numbers. This appears to
be the first direct proof of this fact. The problem of finding such a proof was
suggested by Stanley \cite{9}.
\endabstract
\endtopmatter
\document

\def\mysec#1{\bigskip\centerline{\bf #1}\message{ * }\nopagebreak\par}

\def\myref#1{\item"{[{\bf #1}]}"} 
 
\def\pf{{\it Proof.\ }} 

\def\cite#1{\relaxnext@
  \def\nextiii@##1,##2\end@{[{\bf##1},\,##2]}%
  \in@,{#1}\ifin@\def\next{\nextiii@#1\end@}\else
  \def\next{[{\bf#1}]}\fi\next}
\def\proclaimheadfont@{\smc}

\def\pf{{\it Proof.\ }}
\define\N{{\bold N}}
\define\Z{{\bold Z}}

\define\twoline#1#2{\line{\hfill{\smc #1}\hfill{\smc #2}\hfill}}

\def\mypic#1{\epsffile{#1}}

\mysec{1. Introduction}

A plane partition $\pi$ is a rectangular array of non-negative integers with
non-increasing rows and columns and finitely many nonzero entries. One can naturally
identify $\pi$ with an order ideal of ${\N}^3$, i.e., a finite subset of ${\N}^3$ 
such that $(i,j,k)\in\pi$ implies $(i',j',k')\in\pi$, whenever $i\geq i'$, 
$j\geq j'$ and $k\geq k'$.

By permuting the coordinate axes, one obtains an action of $S_{3}$ on the set
of plane partitions. Let $\pi\mapsto{\pi}^t$ and $\pi\mapsto{\pi}^r$ denote
the symmetries corresponding to interchanging the $x$- and $y$-axes and to
cyclically permuting the coordinate axes, respectively.
For the set of plane partitions $\pi$ contained in the box 
$B(a,b,c):=\{(i,j,k)\in{\N}^3:i<a,j<b,k<c\}$, there is an additional symmetry

$$\pi\mapsto{\pi}^c:=\{(i,j,k)\in{\N}^3:(a-i-1,b-j-1,c-k-1)\notin\pi\},$$
called complementation.

These three symmetries generate a group isomorphic to the dihedral group of
order 12, which has 10 conjugacy classes of subgroups. These lead to 10 
distinct enumeration
problems: determine the number of plane partitions contained in a given box
that are invariant under the action of one of these subgroups. The program of
solving these problems was formulated by Stanley \cite{9} and has been 
recently completed (see \cite{1}, \cite{6} and \cite{11}). 

Even so, there are still many aspects of this group of enumeration questions that 
continue to attract a lot of attention. Although significant progress has been made in
giving unified proofs of some of the cases (see for example \cite{7},\cite{12}), 
there is still no good explanation as to why all cases are enumerated by simple 
product formulas. 

One possible attempt to explain this is to prove simple relations between the numbers
enumerating these classes, without explicitly evaluating them
(see \cite{3} for an illustration of this).

Two of the cases that turned out to be among the hardest to prove are those 
of cyclically 
symmetric, self-complementary plane partitions (i.e., plane partitions $\pi$ with 
${\pi}^r={\pi}^c={\pi}$), first proved by Kuperberg \cite{6}, and totally symmetric, 
self-complementary plane partitions 
(i.e., plane partitions invariant under the full symmetry group of the box), first
proved by Andrews \cite{1}. It is 
easy to see that in order for plane partitions in either of these symmetry classes
to exist, the box must be a cube of even side. Denote by $CSSC(2n)$ and $TSSC(2n)$ the
number of plane partitions in the two classes, respectively, where $2n$ is the 
side of the cubical box.  

In this paper we prove that the former of these two numbers equals the square of the
latter, without explicitly evaluating either of them. The problem of finding such a 
direct proof was suggested by Stanley \cite{9}. One can view our result as providing
new proofs for the two symmetry classes it relates. 

\mysec{2. Proof of the result} 

\proclaim{Theorem 2.1} $CSSC(2n)=TSSC(2n)^2$.
\endproclaim

Our proof employs the following preliminary result.

\proclaim{Proposition 2.2} Let $U(n)$ be the matrix 

$$U(n)=\left(\frac{1}{2}{i+j \choose 2i-j}+{i+j \choose 2i-j-1}\right)_{0\leq i,j\leq
n-1}.\tag2.1$$
Then we have

$$CSSC(2n)=2^n\det U(n).\tag2.2$$
\endproclaim

\pf Consider the tiling of the plane by unit equilateral triangles.
Define a {\it region} to be the union of finitely many such
unit triangles. Let $H(a,b,c)$ be the hexagonal region having sides of lengths
$a,b,c,a,b,c$ (in cyclic order). Then it is well-known (see \cite{4}, \cite{6} 
and \cite{8}) that plane partitions fitting inside $B(a,b,c)$ can be identified with
tilings of $H(a,b,c)$ by unit rhombi (also called {\it lozenge tilings}). Moreover,
all symmetry classes of plane partitions get identified with classes of
lozenge tilings invariant under certain symmetries of $H(a,b,c)$. In particular, 
$CSSC(2n)$ turns out to be equal to the number of lozenge tilings of
$H_n:=H(2n,2n,2n)$ that are invariant under rotation by 60 degrees (see \cite{6}).

In turn, lozenge tilings of $H_n$ can be identified with perfect matchings 
of the dual graph $G_{n}$, i.e., the graph whose vertices are the unit triangles 
contained in $H_n$, and whose edges connect precisely those unit triangles that 
share an edge (a perfect matching of a graph is a collection of vertex-disjoint edges
collectively incident to all vertices of the graph; we usually refer to a perfect
matching simply as a {\it matching}). We obtain that $CSSC(2n)$ equals the number of 
matchings of $G_{n}$ invariant under the rotation $\rho$ by 60 degrees around the 
center of $G_{n}$.

Consider the action of the group generated by $\rho$ on $G_{n}$, and let
$\tilde{G}_{n}$ be the orbit graph. It is easy to see that the 60 degree
invariant matchings of $G_{n}$ can be identified with the matchings of
$\tilde{G}_{n}$. 

\topinsert
\twoline{\mypic{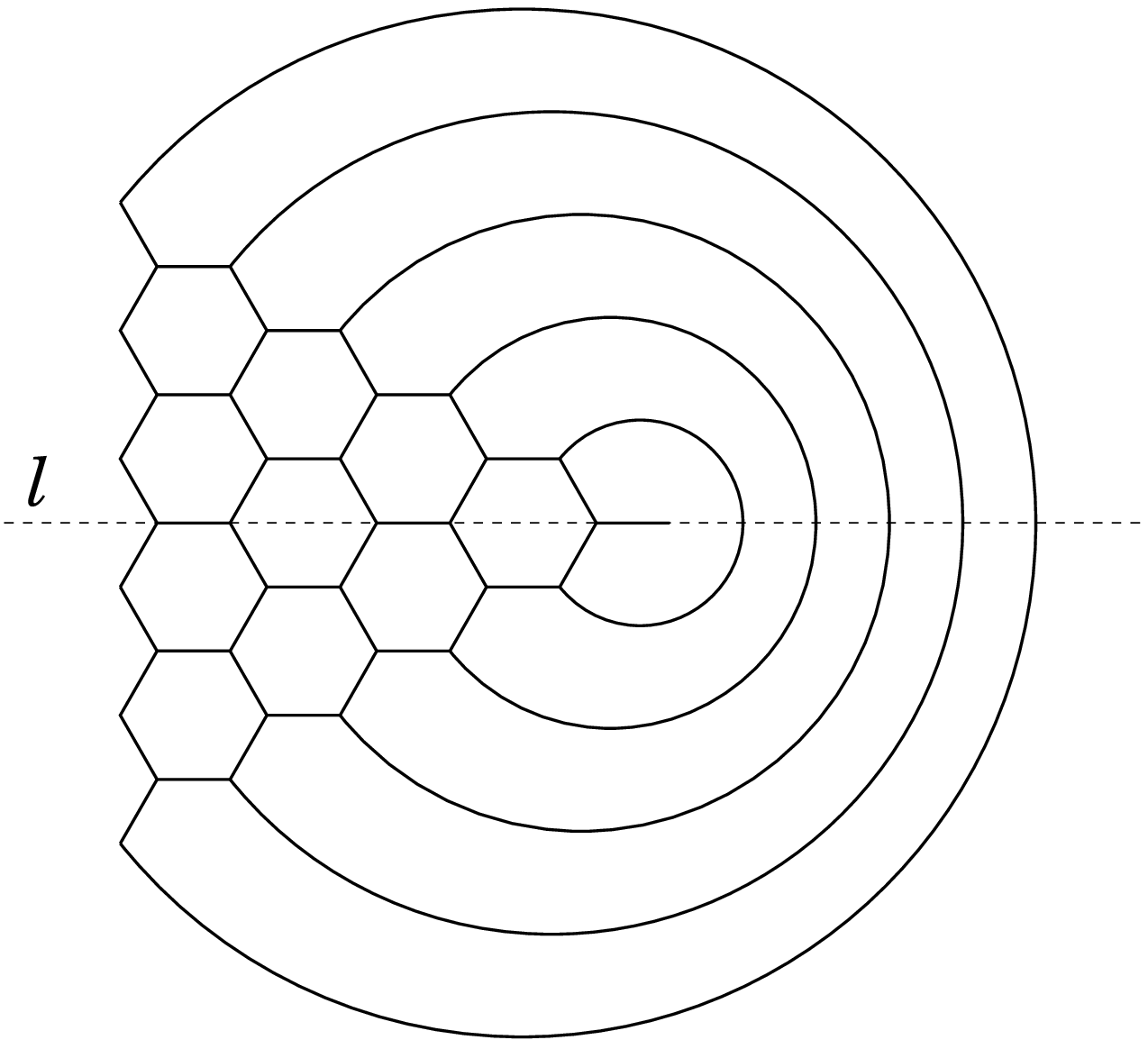}}{\mypic{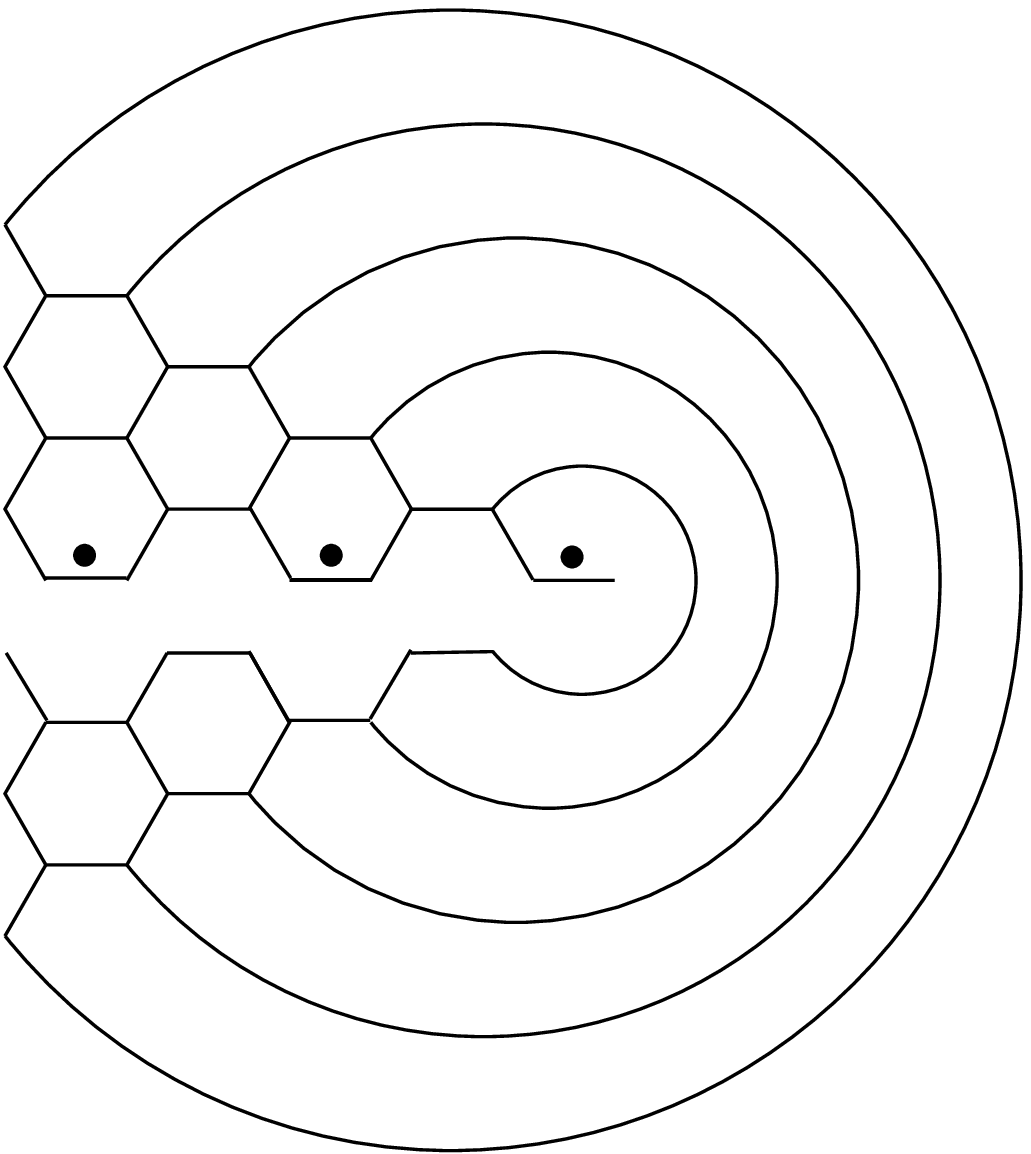}}
\twoline{Figure~2.1. {\rm $\tilde{G}_3$.}}{Figure~2.2. {\rm $K_3$.}}
\endinsert

\topinsert
\twoline{\mypic{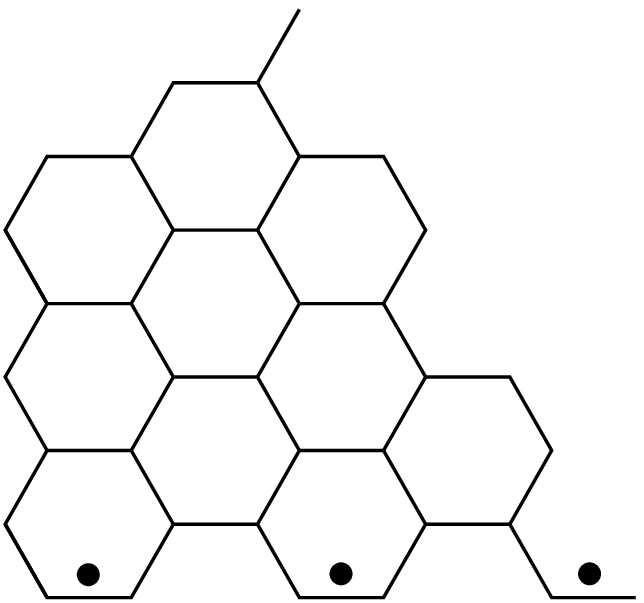}}{\mypic{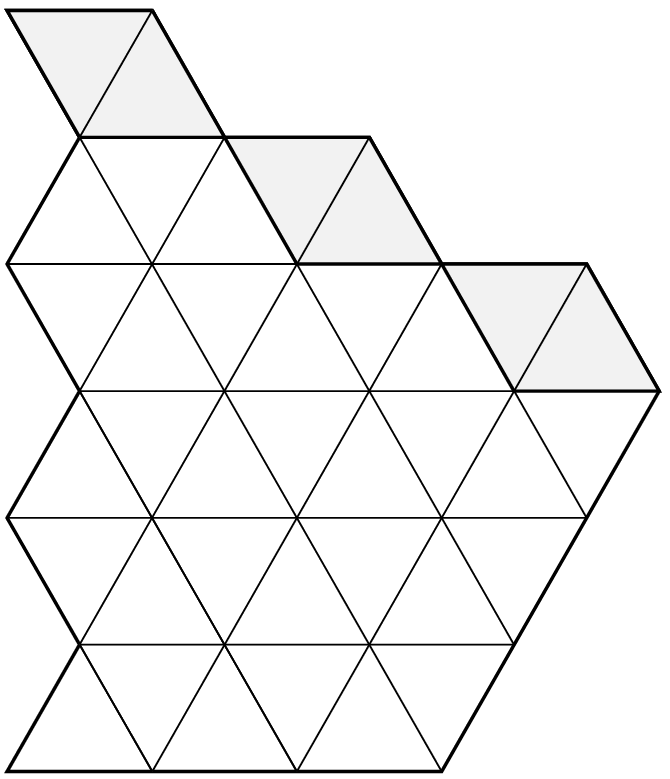}}
\twoline{Figure~2.3. {\rm Another drawing of $K_3$.}}{Figure~2.4. {\rm $R_3$.}}
\endinsert

As shown in Figure 2.1 (for $n=3$), the graph $\tilde{G}_{n}$ 
can be embedded in the plane so
that it admits a symmetry axis $\ell$. (Strictly speaking, $\tilde{G}_{n}$ is the graph
obtained from the one shown in Figure 2.1 by adding a loop to the vertex of degree
one; however, since $\tilde{G}_{n}$ contains no loop besides this and it has an
even number of vertices, this loop is not part of any perfect matching, so it can be
ignored). 

It can be easily checked that the variant of the
Factorization Theorem \cite {3,Theorem 1.2} for matchings presented in 
\cite{3,Proof of Theorem 7.1} can be
applied to $\tilde{G}_{n}$. One obtains that the number of matchings of
$\tilde{G}_{n}$ equals $2^n$ times the matching generating function of the 
subgraph $K_n$ (illustrated in Figure 2.2, for $n=3$) 
obtained by deleting the $2n-1$ edges immediately below $\ell$, 
and changing the weight of the $n$ edges along $\ell$ to 1/2 (the matching generating
function of a graph is the sum of the weights of all its perfect matchings, where the
weight of a matching is the product of weights of its edges). 

The graph $K_n$ can be clearly redrawn in the plane as shown in Figure 2.3. Using  
again the duality between matchings and lozenge tilings, the matchings of $K_n$ can be
identified with tilings of the dual region $R_n$ shown (for $n=3$) in Figure 2.4
(indeed, the dual graph of $R_n$ is the same as the image of $K_n$ under
counterclockwise rotation by 150 degrees). Consider the $n$ tile positions in $R_n$
along its northeastern boundary (they are indicated by a shading in Figure 2.4). In a
tiling of $R_n$, weight each tile occupying one of these positions by 1/2, and all
others by 1; let $L^*(R_n)$ be the tiling generating function of $R_n$ under this
weighting. With this convention, the bijection between matchings of $K_n$ and
tilings of $R_n$ is weight-preserving. Therefore, one obtains

$$CSSC(2n)=2^nL^*(R_n).\tag2.3$$

\topinsert
\twoline{\mypic{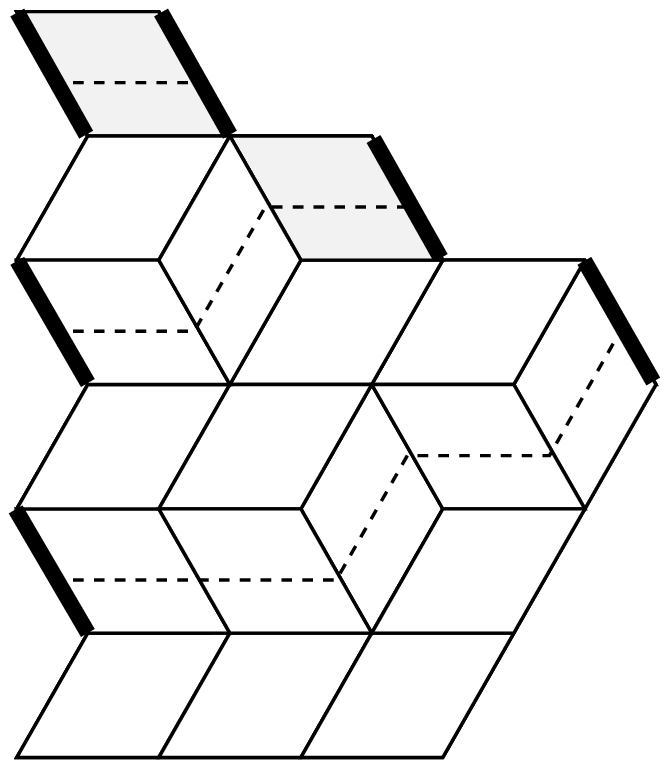}}{\mypic{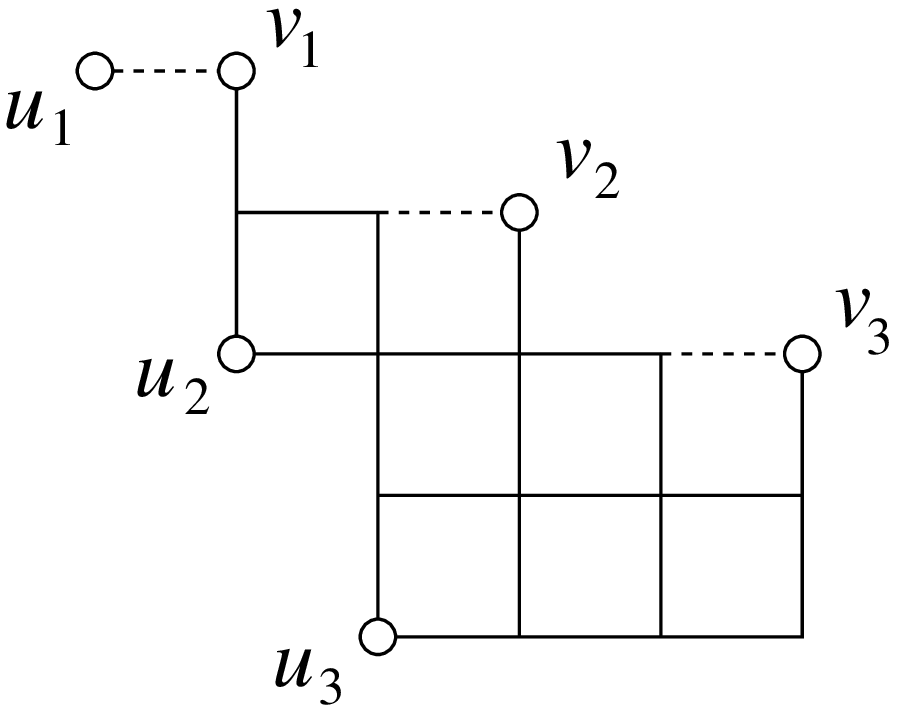}}
\twoline{Figure~2.5}{Figure~2.6}
\endinsert

Let $\Cal T$ be a tiling of $R_n$ and consider the $n$ unit segments facing northeast 
on the left boundary of $R_n$ (these are outlined in thick solid lines in Figure
2.5). By ``following'' the lozenges of $\Cal T$ containing
these segments, one obtains an $n$-tuple of non-intersecting paths of rhombi
connecting our $n$ unit segments to the $n$ unit segments facing southwest on the
right boundary of $R_n$ (also shown in thick solid lines in Figure 2.5; the paths
of rhombi are indicated by dashed lines). This in turn is readily identified with an 
$n$-tuple $\Cal P=(P_0,\dotsc,P_{n-1})$ of non-intersecting paths on the $\Z^2$ 
lattice, where $P_i$ runs from
$u_i=(i,2n-2i)$ to $v_i=(2i+1,2n-i)$, taking unit steps north and east, for
$i=0,\dotsc,n-1$ (see Figure 2.6). Moreover, it is not hard to see that $\Cal P$ 
determines $\Cal T$ uniquely.

Regard $\Z^2$ as a directed graph, with the edges oriented from west to east
and from south to north. Assign weight 1/2 to the horizontal edge whose right vertex 
is $v_i$, for
$i=0,\dotsc,n-1$, and weight 1 to all other edges of $\Z^2$ (the edges weighted by 1/2
are showed in dotted lines in Figure 2.6). Define the weight 
of a lattice path to be the product of the weights of its steps. The weight of a
$k$-tuple of lattice paths is defined as the product of the weights of the individual
paths. The generating function of a family of $k$-tuples of paths is the sum of
weights of its members.

By our choice of weights on $\Z^2$ and by the above-described bijection between 
tilings and non-intersecting lattice paths, it is clear that  $L^*(R_n)$ is the 
generating function of $n$-tuples of non-intersecting lattice paths 
$\Cal P=(P_0,\dotsc,P_{n-1})$, where $P_i$ runs from $u_i$ to $v_i$ ($i=0,\dotsc,n-1$).

Since our orientation of $\Z^2$ is acyclic and since the $n$-tuples $(u_i)$ and
$(v_i)$ of starting and ending points of our paths are compatible in the sense of
\cite{10,Theorem1.2}, one obtains by Theorem 1.2 of \cite{10} (see also \cite{5})
that the generating
function for our $n$-tuples $\Cal P$ of non-intersecting lattice paths equals the
determinant of the $n\times n$ matrix whose $(i,j)$-entry is the generating 
function of lattice paths from $u_i$ to $v_j$, $i,j=0,\dotsc,n-1$. From our choice of
weights, it is readily seen that this is precisely the $(i,j)$-entry of the matrix 
$U(n)$ given by (2.1). 
It follows that $L^*(R_n)=\det U(n)$, and hence using (2.3) we obtain (2.2) $\square$

\smallpagebreak
{\it Proof of Theorem 2.1.} We use a result of Stembridge \cite{10} which expresses
$TSSC(2n)$ as a Pfaffian of order $n$. In a restatement due to Andrews \cite{1}
this result is

$$TSSC(2n)^2=\det st(n),\tag2.4$$
where the matrix $st(n)=(a_{ij})_{0\leq i,j\leq n-1}$ is given by

$$
a_{ij}=\left\{
\aligned &1\ \ \ \ \ \ \ \ \ \ \ \ \ \ \ \ \ \ \ \ \ \ \ \ \ \text{if}\ i=j=0,\\ 
&0\ \ \ \ \ \ \ \ \ \ \ \ \ \ \ \ \ \ \ \ \ \ \ \ \ \text{if}\  i=j>0,\\
&\sum_{s=2i-j+1}^{2j-i}{i+j\choose s}\ \ \ \ \text{if}\  i<j,\\
&-a_{ji}\ \ \ \ \ \ \ \ \ \ \ \ \ \ \ \ \ \ \ \text{if}\  i>j.
\endaligned
\right.
$$



In \cite{2} Andrews and Burge show (see relations (4.12) and (4.13) of \cite{2}), 
by means of simple row and column operations, that 

$$\det st(n)=\det w(n),\tag2.5$$ 
where

$$w(n)=\left({i+j+1 \choose 2i-j}+{i+j \choose 2i-j-1}\right)_{0\leq i,j\leq
n-1}.\tag2.6$$
It follows from (2.4) and (2.5) that

$$TSSC(2n)^2=\det w(n).\tag2.7$$
However, comparing (2.1) and (2.6), one readily checks that each entry of $U(n)$ is
precisely half of the corresponding entry of $w(n)$. Therefore, by (2.2) and (2.7) we
obtain $CSSC(2n)=TSSC(2n)^2$, which completes the proof. $\square$

\mysec{References}
{\openup 1\jot \frenchspacing\raggedbottom
\roster
\myref{1}
  G. E. Andrews, Plane Partitions, V: The T.S.S.C.P.P. conjecture, {\it J.
Comb. Theory Ser. A} {\bf 66} (1994), 28--39.
\myref{2}
  G. E. Andrews and W. H. Burge, Determinant identities, {\it Pacific J. Math.}, {\bf
158} (1993), 1--14.
\myref{3}
  M. Ciucu, Enumeration of perfect matchings in graphs with reflective symmetry, 
{\it J. Comb. Theory Ser. A} {\bf 77} (1997), 67--97.
\myref{4}
  G. David and C. Tomei, The problem of the calissons, {\it Amer. Math. 
Monthly} {\bf 96} (1989), 429--431.
\myref{5}
  I. M. Gessel and X. Viennot, Determinants, paths, and plane partitions, preprint,
1989.
\myref{6}
  G. Kuperberg, Symmetries of plane partitions and the permanent-determinant
method, {\it J. Comb. Theory Ser. A} {\bf 68} (1994), 115--151.
\myref{7}
  G. Kuperberg, Four symmetry classes of plane partitions under one roof, {\it J. 
Combin. Theory Ser. A} {\bf 75} (1996), 295-315.
\myref{8}
  D. P. Robbins, The story of $1,2,7,42,429,7436,\dotsc$, {\it Math. Intelligencer}
{\bf 13} (1991), 12--19.
\myref{9}
  R. P. Stanley, Symmetries of plane partitions, {\it J. Comb. Theory Ser. A} 
{\bf 43} (1986), 103--113.
\myref{10}
  J. R. Stembridge, Nonintersecting paths, Pfaffians and plane partitions, 
{\it Adv. in Math.} {\bf 83} (1990), 96--131.
\myref{11}
  J. R. Stembridge, The enumeration of totally symmetric plane partitions,
{\it Adv. in Math.} {\bf 111} (1995), 227--243.
\myref{12}
  J. R. Stembridge, Some hidden relations involving the ten symmetry classes of plane
partitions, {\it J. Comb. Theory Ser. A} {\bf 68} (1994), 372--409.

\endroster\par}

\enddocument